\input amstex
\input amssym
\documentstyle{amsppt}
\magnification=1200
\magnification=\magstep1
\hsize=6truein
\vsize=8.4truein
\baselineskip 0.7cm
\parindent=0.8cm

 \NoBlackBoxes \TagsOnLeft \magnification=\magstep1
\centerline {\bf Action of Dihedral Groups}

 \centerline {Ming-chang Kang}

\centerline {Department of Mathematics}

\centerline {National Taiwan University \,}

\centerline {  Taipei, Taiwan \qquad \qquad \qquad }

\centerline {\quad E-mail: kang\@math.ntu.edu.tw}

\

\document
\baselineskip 0.7cm
\parindent=0.8cm

\

\noindent {\bf Abstract. }  Let $K$ be any field and $G$ be a finite group.  Let $G$ act on the
rational function field $K(x_g: \ g \in G)$ by $K$-automorphisms defined by $g \cdot x_h=x_{gh}$
for any $g, \ h \in G$.  Denote by $K(G)$ the fixed field $K(x_g: \ g \in G)^G$.  Noether's problem
asks whether $K(G)$ is rational (=purely transcendental) over $K$.  We will give a brief survey of
Noether's problem for abelian groups and dihedral groups, and will show that $\Bbb Q(D_n)$ is
rational over $\Bbb Q$ for $n \le 10$.

\vskip 2cm

\noindent $\underline {\mskip 200mu}$

\noindent Mathematics Subject Classification 2000: Primary 12F12, 13A50, 11R32, 14E08.

\noindent Keywords: Noether's problem, rationality problem, dihedral groups.

\noindent This article was written for a local conference in 2005.
It was circulated among a few friends, but has not been published
ever since. It is not difficult to adapt the proof of this article
so that the base field $\Bbb Q$ is replaced by a rather general
field $k$.

\newpage

\noindent {\bf \S 1. Introduction}

Let $K$ be any field and $G$ be a finite group.  Let $G$ act on the rational function field $K(x_g:
\ g \in G)$ by $K$-automorphisms such that $g \cdot x_h=x_{gh}$ for any $g, \ h \in G$.  Denote by
$K(G)$ the fixed field $K(x_g: \ g \in G)^G=\{f \in K(x_g: g \in G): \ \sigma \cdot f = f \ {\text
{for any}} \ \sigma \in G\}$.  Noether's problem asks whether $K(G)$ is rational (=purely
transcendental) over $K$ \cite {No}.

Noether's problem arose from the study of the inverse Galois problem.  In particular,
if $K(G)$ is rational and $K$ is an infinite field, then a generic polynomial for
Galois $G$-extensions over $K$ exists \cite {DM; Sa1}. In case $K$ is a Hilbertian
field, i.e. Hilbert irreducibility theorem is valid for irreducible polynomials $f \in
K[x_1, x_2, \cdots, x_n]$ (for example, an algebraic number field or a field $K$ which
is finitely generated over some field $k$ so that trans deg$_k K \ge 1$), the existence
of a generic polynomial for Galois $G$-extensions over $K$ will certainly guarantee an
infinite family of Galois field extensions of $K$ with Galois groups isomorphic to $G$.

The first solution of Noether's problem is provided by E. Fischer, a friend of Emmy
Noether introducing to her the then novel and abstract thinking of Dedekind and
Hilbert.

\proclaim {Theorem 1.1} {\rm {(Fischer \cite {Fi})}} Let $G$ be a finite abelian group with
exponent $e$ and $K$ be any field containing a primitive $e$-th root of unity.  Then $K(G)$
 is rational over $K$. \endproclaim

On the other hand only a handful of results for the rationality of $\Bbb Q(G)$ were
known before 1950's.  Samson Breuer was able to show that $\Bbb Q(\Bbb Z_3)$ and $\Bbb
Q(\Bbb Z_6)$ are rational over $\Bbb Q$ where $\Bbb Z_n$ is the cyclic group of order
$n$ \cite {Br1}; he then showed that $\Bbb Q(G)$ is rational for some transitive
solvable subgroup $G$ contained in $S_p$, the symmetric group of degree $p$, if $p=5$
or 7 \cite {Br2}.  Breuer's results for transitive solvable subgroups in $S_p$ was
extended by Furtw\"angler for $p=5, \ 7, \ 11$ \cite {Fr}; finally Breuer himself
extended these results for any prime number $p \le 23$ \cite {Br3}. Several years later
Gr\"obner proved that $\Bbb Q(G)$ is rational if $G$ is the quaternion group of order 8
\cite {Gr}.

Unfortunately almost all these results, except Fischer's Theorem, were forgotten after
World War II.  In 1955 H. Kuniyoshi and K. Masuda resumed this problem ; they called it
a problem of Chevalley. Masuda rediscovered many previous results; in particular he
proved that $\Bbb Q(\Bbb Z_p)$ is rational if $p=3, \ 5, \ 7, \ 11$ \cite {Ma}. To the
surprise of most people Swan constructed the first counter-example to Nother's problem
in 1969 \cite {Sw1}: $\Bbb Q(\Bbb Z_p)$ is not rational over $\Bbb Q$ if $p=47, \ 113,
\ 233, \cdots$.  The reader is referred to the survey articles \cite {Sw1; Sw2; Sa2;
Ke; Ka2} for subsequent progress of Noether's problem.  We will remark that Noether's
problem for finite abelian groups has been solved completely \cite {Le}.

In the remaining of this article we will focus on the rationality of $\Bbb Q(D_n)$
where $D_n= < \sigma, \ \tau: \ \sigma ^n=\tau^2=1, \ \tau \sigma \tau^{-1}=\sigma
^{-1}>$ is the dihedral group of order $2n$.  As mentioned before, the rationality of
$\Bbb Q(D_p)$ was proved for prime numbers $p \le 23$ by S. Breuer and Furtw\"angler.
It seems strange that the answer to the rationality of $\Bbb Q(D_n)$, $3 \le n \le 10$,
is difficult to locate in the literature.  The rationality of $\Bbb Q(D_8)$ can be
found in \cite {CHK, Theorem 3.1}; an easier case, the rationality of $\Bbb Q(D_4)$, is
given in \cite {CHK, Proposition 2.6}.  In fact, $K(G)$ is rational for any field $K$
and any non-abelian group $G$ of order 8 or 16, except for the case when $G$ is the
generalized quaternion group of order 16 \cite {CHK; Ka2}.  The task of this article is
to study the rationality problem of $\Bbb Q(D_6)$, $\Bbb Q(D_9)$ and $\Bbb Q(D_{10})$.
What we will prove is the following.

\proclaim {Theorem 1.2}  $\Bbb Q(D_n)$ is rational over $\Bbb Q$ for $3
\le n \le 10$. \endproclaim

Needless to say, the method in proving the above theorem may be adapted for a more
general context, which will be embodied in a separate article. One of the purposes of
this article is to illustrate some techniques of solving Noether's problem through
concrete cases.

The article is organized as follows:  Some basic tools will be recalled in
Section 2.  The rationality of $\Bbb Q(D_9)$ (resp. $\Bbb Q(D_6)$ and
$\Bbb Q(D_{10})$) will be proved in Section 3 (resp. Section 4).  The
proof of Theorem 1.2 will be finished once we obtain Theorem 3.1, Theorem
4.1 and Theorem 4.2.

\

\noindent {\bf \S 2.  Preliminaries}

We recall some basic results which will be used in Section 3 and Section
4.

\proclaim {Theorem 2.1} {\rm {([CHK, Theorem 2.1])}}  Let $L$ be a field
and $G$ be a finite group acting on $L(x_1, \cdots, x_n)$, the rational
function field of $m$ variables over $L$.  Suppose that \roster \item
"(i)" for any $\sigma \in G$, $\sigma (L) \subset L$;
\item "(ii)" the restriction of the action of $G$ on $L$ is faithful;
\item "(iii)" for any $\sigma \in G$, $$\pmatrix \sigma(x_1)\\ \vdots \\
\sigma (x_m) \endpmatrix= A(\sigma)\pmatrix x_1\\ \vdots \\ x_m
\endpmatrix+B(\sigma)$$\endroster
where $A(\sigma) \in GL_m(L)$ and $B(\sigma)$ is an $m \times 1$ matrix over $L$.  Then
there exists $z_1, \cdots, z_m \in L(x_1, \cdots, x_m)$ with $L(x_1, \cdots, x_m) =
L(z_1, \cdots, z_m)$ such that $\sigma(z_i)=z_i$ for any $\sigma \in G$, any $1 \le i
\le m$.\endproclaim

\proclaim {Theorem 2.2} {\rm {([CHK, Theorem 2.4])}} Let $G$ be any group whose order
may be finite or infinite. Suppose that $G$ acts on $L(x)$, the rational function field
of one variable over a field $L$. Assume that, for any $\sigma \in G, \, \sigma(L)
\subset L$, and $\sigma(x) = a_{\sigma} \cdot x + b_{\sigma}$ for some $a_{\sigma},
b_{\sigma} \in L$ with $a_{\sigma} \neq 0$. Then $L(x)^G$ is rational over $L^G$.
\endproclaim

\proclaim {Theorem 2.3} {\rm {([CHK, Theorem 2.3])}}  Let $K$ be any field, $a, \ b \in
K \setminus \{0\}$ and $\sigma: K(x,y) \longrightarrow K(x,y)$ be a $K$-automorphism
defined by $\sigma(x)=a/x$, $\sigma(y)=b/y$. Then $K(x, y)^{<\sigma>}=K(u,v)$ where
$$\displaystyle u=\frac {x-\dfrac ax}{xy-\dfrac {ab}{xy}}, \qquad v=\frac {y-\dfrac
by}{xy-\dfrac {ab}{xy}}.$$ \endproclaim

\proclaim {Theorem 2.4}{\rm {(Hajja [Ha])}}  Let $K$ be any field and $G$ be a finite
group.  Suppose that $G$ acts on the rational function field $K(x_1, x_2)$ by
$K$-automorphisms such that, for any $\sigma \in G$, for $1 \le i \le 2$,
$\sigma(x_i)=a_i(\sigma)\cdot x_1^{m_i(\sigma)}x_2^{n_i(\sigma)}$ where $a_i(\sigma)
\in K \setminus \{0\}$ and $m_i(\sigma), \ n_i(\sigma) \in \Bbb Z$.  Then $K(x_1,
x_2)^G$ is rational over $K$. \endproclaim

\

\noindent {\bf \S 3. $\Bbb Q(D_9)$}

Let $D_n=<\sigma, \ \tau: \sigma^n=\tau^2=1, \tau \sigma
\tau^{-1}=\sigma^{-1}>$ be the dihedral group of order $2n$.  Let $V=
\bigoplus _{g \in D_n} \Bbb Q \cdot x(g)$ be the regular representation
space of $D_n$ over $\Bbb Q$, i.e. $g \cdot x(h)=x(gh)$ for any $g, \ h
\in D_n$.

Define $x_i= x(\sigma^i)+x(\sigma^i \tau)$ for $0 \le i \le n-1$.  Then
$$\aligned \sigma&: x_0 \mapsto x_1 \mapsto \cdots \mapsto x_{n-1} \mapsto
x_0, \\ \tau &: x_i \mapsto x_{-i} \endaligned $$ where the index of $x_i$
is taken modulo $n$.

Clearly $\bigoplus _{0 \le i \le n-1} \Bbb Q \cdot x_i$ is a faithful
$D_n$-subspace of $V$.  By Theorem 2.1, if $\Bbb Q(x_0, \cdots,
x_{n-1})^{D_n}$ is rational over $\Bbb Q$, then $\Bbb Q(D_n)=\Bbb Q(x_g: g
\in D_n)^{D_n}$ is also rational over $\Bbb Q$.

We will prove that $\Bbb Q(x_0, \cdots, x_{n-1})^{D_n}$ is rational over
$\Bbb Q$ for $n=9, 6, 10$.

\proclaim {Theorem 3.1}  $\Bbb Q(x_0, \cdots, x_{8})^{D_9}$ is rational
over $\Bbb Q$. \endproclaim

\demo {Proof} Step 1.  Let $\zeta$ be a primitive 9-th root of unity and
$\pi=Gal(\Bbb Q(\zeta)/\Bbb Q)$. Let $\rho \in \pi$ such that
$\rho(\zeta)=\zeta^2$. Then $\rho$ is a generator of $\pi$.

Extend the actions of $D_9$ and $\pi$ to $\Bbb Q(\zeta)(x_0, \cdots,x_8)$ by
requiring $\sigma(\zeta)=\tau(\zeta)=\zeta$, $\rho(x_i)=x_i$ for $0 \le i \le
8$.

Define $y_i= \sum _{0 \le j\le8} \zeta^{-ij}x_j$.  Then $$\aligned \sigma &:
y_i\mapsto \zeta^iy_i, \\ \tau &:y_i\mapsto y_{-i}, \\ \rho&: y_i \mapsto
y_{2i}
\endaligned$$

We find that $\Bbb Q(x_0, \cdots, x_8)^{D_9}=\{\Bbb Q(\zeta)(x_0, \cdots,
x_8)^{<\rho>}\}^{D_9}=\Bbb Q(\zeta)(x_0, \cdots$, \, $x_8)^{<D_9, \rho>}$ $=\Bbb
Q(\zeta)(y_0, \cdots, y_8)^{<\sigma, \tau, \rho>}$.  Moreover $\tau \cdot
\rho^3(y_i)=y_i$ for any $0 \le i \le 8$.

By Theorem 2.1, if $\Bbb Q(\zeta)(y_i: i \in \Bbb Z_9^\times)^{<\sigma, \tau,
\rho>}$ is rational over $\Bbb Q$, so is $\Bbb Q(\zeta)(y_0, \cdots,$ $
y_8)^{<\sigma, \tau, \rho>}$.  Hence it suffices to prove that  $\Bbb
Q(\zeta)(y_i: i \in \Bbb Z_9^\times)^{<\sigma, \tau, \rho>}$ is rational over
$\Bbb Q$.

Step 2.  Let $<y_1, y_2, y_4, y_5, y_7, y_8>$ be the multiplicative subgroup of $\Bbb
Q(\zeta)(y_0$, $\cdots, y_8)\setminus \{0\}$ generated by $y_i$ where $i \in \Bbb
Z_9^{\times}$.  As an abelian group, the group $<y_1, y_2, y_4$, $y_5, y_7, y_8>$ is
isomorphic to the free abelian group $\Bbb Z^6$.

There is a natural $\pi$-module structure on $<y_1, y_2, y_4, y_5, y_7, y_8>$. As a
module over $\Lambda:=\Bbb Z[\pi]$, $<y_1, y_2, y_4, y_5, y_7, y_8>$ is isomorphic to
$\Lambda$.  In fact,we may identify $y_1, y_2$ with $1, \rho \in \Lambda$ respectively
and write $\Lambda = <y_1, y_2, y_4, y_5, y_7, y_8>$.

Define a map $$\aligned \Phi: \Lambda=<y_1, y_2, y_4, y_5, y_7,
y_8&>\longrightarrow \Bbb Z_9\\ & y \longmapsto \overline  {j}\endaligned$$ if
$\sigma(y)=\zeta^jy$ where
$y=y_1^{a_1}y_2^{a_2}y_4^{a_4}y_5^{a_5}y_7^{a_7}y_8^{a_8}$ with $a_i \in \Bbb
Z$.

Define a $\pi$-module structure on $\Bbb Z_9$ by $\rho \cdot \overline j =
\overline {2j}$.  Thus $\Phi$ becomes a $\pi$-equivariant map, i.e.
$\Phi(\lambda \cdot y)=\lambda \cdot \Phi(y)$ for any $\lambda \in \pi$.

Let $M$ be the kernel of $\Phi$. Then $\Bbb Q(\zeta)(y_1, y_2, y_4, y_5, y_7,
y_8)^{<\sigma>}=\Bbb Q(\zeta)(M)$. Note that $M$ is an ideal of $\Lambda$ with
$[\Lambda :M]=9$.

By \cite {Le, (3.3) Proposition, p. 311}, $M$ is a projective $\Lambda$-module.

Step 3.  We will prove that $M$ is a free module over $\Lambda$.  In fact, we
will show that any rank-one $\Lambda$-projective module is free.

Write $\Lambda=\Bbb Z[T]/<T^6-1>$, $\Lambda_1=\Bbb Z[T]/<T^3+1>$, $\Lambda_2=\Bbb
Z[T]/<T^3-1>$. (Note that $\Lambda _1 \simeq \Lambda_2$.)  It is not difficult to prove
that all the Picard groups of $\Lambda, \ \Lambda_1, \ \Lambda_2$ are zero by using
Mayer-Vietoris sequences of $K$-groups \cite {Mi, Theorem 3.3, p. 28} for the Cartesian
squares $$\displaystyle \CD \Lambda @> >> \Lambda_1 \\@V VV @VV V
\\ \Lambda_2 @> >> {\Bbb Z_2[T]}/ {<T^3-1>} \endCD \qquad \CD \Lambda_2 @> >> \Bbb
Z[\zeta_3]
\\@V VV @VV V \\ \Bbb Z @> >> \Bbb Z_3 \endCD$$
where $\zeta_3$ is a primitive 3rd root of unity.

Step 4.  By Step 3, we may find $z_0 \in <y_1, y_2, y_4, y_5, y_7, y_8>$ such
that $\Bbb Q(M)= \Bbb Q(z_0, z_1, \cdots, z_5)$where $z_i=\rho^i(z_0)$ for $1
\le i \le 5$.  Clearly $\rho: z_0 \mapsto z_1 \mapsto z_2 \mapsto z_3 \mapsto
z_4 \mapsto z_5 \mapsto z_0$ and $\tau \rho^3$ is the identity map on $\Bbb
Q(z_0, \cdots, z_5)$.

Now $\Bbb Q(\zeta)(y_1, y_2, y_4, y_5, y_7, y_8)^{<\tau \rho^3>}=\Bbb Q(\zeta)(z_0,
\cdots, z_5)^{<\tau \rho^3>}= \Bbb Q(\eta)(z_0, \, \cdots$, $z_5)$ where $\eta= \zeta
+\zeta^{-1}$.  It remains to show that $\Bbb Q(\eta)(z_0, \cdots, z_5)^{<\rho>}$ is
rational.

Define $u_i=z_i-z_{i+3}$, $v_i=z_i+z_{i+3}$ for $0 \le i \le 2$.  We find that
$$\rho: u_0 \mapsto u_1 \mapsto u_2 \mapsto -u_0, \quad v_0 \mapsto v_1 \mapsto
v_2 \mapsto v_0.$$

By Theorem 2.1, it suffices to show that $\Bbb Q(\eta)(u_0, u_1,
u_2)^{<\rho>}$ is rational over $\Bbb Q$.

Step 5.  Note that $\rho^3(u_i)=-u_i$.  Define $v_0=u_0^2, \ v_1=u_1/u_0, \
v_2=u_2/u_1$.  Then $\Bbb Q(\eta)(u_0, u_1, u_2)^{<\rho^3>}= \Bbb Q (\eta)(v_0, v_1,
v_2)$.  Moreover, $$\rho: v_0 \mapsto v_0v_1^2, \quad v_1 \mapsto v_2 \mapsto
-1/(v_1v_2).$$

By Theorem 2.1, if $\Bbb Q(\eta)(v_1, v_2)^{<\rho>}$ is rational over
$\Bbb Q$, then $\Bbb Q(\eta)(v_0, v_1, v_2)^{<\rho>}$ is also rational
over $\Bbb Q$.

Step 6. Define $w_1=1/(1-v_1+v_1v_2), \ w_2=-v_1/(1-v_1+v_1v_2)$.  Then
$\Bbb Q(\eta)(v_1, v_2)=\Bbb Q(\eta)(w_1, w_2)$ and $\rho: w_1 \mapsto w_2
\mapsto -w_1-w_2+1$.  By Theorem 2.1, $\Bbb Q(\eta)(w_1, w_2)=\Bbb
Q(\eta)(X,Y)$ for some  $X, \ Y$ such that  $\rho(X)=X, \ \rho(Y)=Y$.
Hence $\Bbb Q(\eta)(w_1, w_2)^{<\rho>}=\Bbb Q(\eta)(X, Y)^{<\rho>}
=\Bbb Q(\eta)^{<\rho>}(X,Y)=\Bbb Q(X,Y)$\qed \enddemo

\definition {Remark}  We may prove that $\Bbb Q(D_n)$ is rational when
$n=3, \ 5, \ 6, \ 7$ by the same method as above.  For the case when $n = 6$, see
Theorem 4.1 for another proof.\enddefinition

\

\noindent {\bf \S 4. $\Bbb Q(D_6)$ and $\Bbb Q(D_{10})$}

We will prove that $\Bbb Q(x_0, \cdots, x_{n-1})^{D_n}$ is rational for
$n=$ 6, 10 where $\sigma: x_0 \mapsto x_1 \mapsto \cdots \mapsto x_{n-1}
\mapsto x_0, \ \tau : x_i \mapsto x_{-i}$.

Let $n=2m$ where $m$ is an odd integer. (In case $n=$6, 10, $m$ is actually an odd
integer.)  Define $y_i=x_i-x_{m+i}, \ y_i'=x_i+x_{m+i}$ where $0 \le i \le m-1$. We get
$$\aligned \sigma &: y_0 \mapsto y_1 \mapsto \cdots \mapsto y_{m-1} \mapsto -y_0, \ \
y_0' \mapsto y_1' \mapsto \cdots \mapsto y'_{m-1} \mapsto y'_0, \\ \tau &: y_0 \mapsto
y_0, \ \  y'_0 \mapsto y'_0, \ \ y_i \mapsto -y_{-i} , \ \ y'_i \mapsto y'_{-i}
\endaligned$$
where the index of $y_i$ is taken modulo $m$.

By Theorem 2.1, it suffices to prove the rationality of $\Bbb Q(y_0,
\cdots, y_{m-1})^{D_n}$.

\proclaim {Theorem 4.1}  $\Bbb Q(y_0, y_1, y_2)^{D_6}$ is rational over $\Bbb Q$.
\endproclaim

\demo {Proof}  Recall that $\sigma:y_0 \mapsto y_1 \mapsto y_2 \mapsto -y_0, \ \tau :
y_0 \mapsto y_0, \ y_1 \mapsto -y_2, \ y_2 \mapsto -y_1$.

Define $z_1=y_1/y_0, \ z_2=y_2/y_1$.  We find $$\aligned \sigma &:y_0 \mapsto y_0z_1, \
z_1 \mapsto z_2 \mapsto -1/(z_1z_2), \\ \tau &: y_0 \mapsto y_0, \ z_1 \mapsto -z_1z_2,
\ z_2 \mapsto 1/z_2.\endaligned$$

By Theorem 2.2, it suffices to prove that $\Bbb Q(z_1, z_2)$ is rational over $\Bbb Q$.
But $\Bbb Q(z_1, z_2)^{<\sigma, \tau>}$ is rational by Theorem 2.4. \qed \enddemo

\proclaim {Theorem 4.2} $\Bbb Q(y_0, \cdots, y_4)^{D_{10}}$ is rational over $\Bbb Q$.
\endproclaim

\demo {Proof}  Recall that $\sigma : y_0 \mapsto y_1 \mapsto y_2 \mapsto y_3 \mapsto
y_4 \mapsto -y_0, \ \tau: y_0 \mapsto y_0, \ y_1 \mapsto -y_4, \ y_2 \mapsto -y_3, \
y_3 \mapsto -y_2, \ y_4 \mapsto -y_1$.

Step 1.  Let $\zeta$ be a primitive 5-th root of unity and $\pi = Gal(\Bbb
Q(\zeta)/\Bbb Q)$.  Let $\rho \in \pi$ such that $\rho(\zeta)=\zeta^2$.
Then $\rho$ is a generator of $\pi$.

Extend the actions of $D_{10}$ and $\pi$ to $\Bbb Q(\zeta)(y_0, \cdots,
y_4)$ by requiring $\sigma(\zeta)=\tau(\zeta)=\zeta, \ \rho(y_i)=y_i$ for
$0 \le i \le 4$.  Note that $\Bbb Q(y_0, \cdots, y_4)^{<\sigma,
\tau>}=\Bbb Q(\zeta)(y_0,$ $\cdots, y_4)^{<\sigma, \tau, \rho>}$.

For $0 \le i \le 4$, define
$z_i=\prod_{j \ne i} (\sigma +\zeta^j)(y_0)$.  We find that $$\aligned
\sigma &: z_i \mapsto -\zeta^iz_i, \\ \tau &: z_i \mapsto
\zeta^{-2i}z_{-i}, \\ \rho &: \zeta \mapsto \zeta^2, \ z_i \mapsto z_{2i}
\endaligned$$
where the index of $z_i$ is taken modulo 5.

By Theorem 2.1, it suffices to prove that $\Bbb Q(\zeta)(z_1, z_2, z_3,
z_4)^{<\sigma, \tau, \rho>}$ is rational over $\Bbb Q$.

Step 2.  For $2 \le i \le 4$, define $u_i=z_i/z_{i-1}$.  Then $\Bbb
Q(\zeta)(z_1, z_2, z_3, z_4)=\Bbb Q(\zeta)(z_1, u_2, u_3, u_4)$ and
$$\aligned \sigma &: z_1 \mapsto -\zeta z_1, \ u_i \mapsto \zeta u_i \
{\text { \ \ for }} 2 \le i \le 4, \\ \tau &: z_1 \mapsto \zeta^3z_1u_2u_3u_4, \
u_2 \mapsto \zeta^3/u_4, \ u_3 \mapsto \zeta^3/u_3, \ u_4 \mapsto
\zeta^3/u_2, \\ \rho &: z_1 \mapsto z_1 u_2, \ u_2 \mapsto u_3 u_4, \ u_3
\mapsto 1/(u_2u_3u_4), \ u_4 \mapsto u_2 u_3. \endaligned $$

By Theorem 2.2, it suffices to prove that $\Bbb Q(\zeta)(u_2, u_3, u_4)^{<\sigma, \tau,
\rho>}$ is rational over $\Bbb Q$.

Step 3.  Define $v_1=u_2^5, \ v_2=u_4/u_2, \ v_3=u_3/u_2$.  Then $\Bbb
Q(\zeta)(u_2, u_3, u_4)^{<\sigma>}=\Bbb Q(\zeta)(v_1, v_2, v_3)$ and
$$\aligned \tau &: v_1 \mapsto 1/(v_1v_2^5), \ v_2 \mapsto v_2, \ v_3
\mapsto v_2/v_3, \\ \rho &: v_1 \mapsto v_1^2v_2^5v_3^5, \ v_2 \mapsto
1/v_2, \ v_3 \mapsto 1/(v_1v_2^2v_3^2). \endaligned$$

Clearly $\Bbb Q(\zeta)(v_1, v_2, v_3)^{<\tau \rho^2>}=\Bbb Q(\eta)(v_1,
v_2, v_3)$ where $\eta=\zeta+\zeta^{-1}$.  Note that $\Bbb Q(\eta)=\Bbb
Q(\sqrt 5)$.  It remains to find $\Bbb Q(\sqrt 5)(v_1, v_2,
v_3)^{<\rho>}$.

Step 4.  Define $t=1/v_2$, $x=v_1v_2v_3^2, \ y=v_3$.  Then $\Bbb Q(\sqrt 5)(v_1, v_2,
v_3)=\Bbb Q(\sqrt 5)(t,x,y)$ and $$\aligned \rho &: \sqrt 5 \mapsto -\sqrt 5, \
t \mapsto 1/t, \ x \mapsto y \mapsto t/x \mapsto 1/(ty) \mapsto x , \\
\rho^2 &: t \mapsto t, \ x \mapsto t/x, \ y \mapsto 1/(ty). \endaligned$$

Define $$u=\frac {x-\dfrac ax}{xy-\dfrac{ab}{xy}}, \qquad v=\frac {y-\dfrac
by}{xy-\dfrac{ab}{xy}}$$ where $a=t, \ b=1/t$.  Apply Theorem 2.3.  We get $\Bbb
Q(\sqrt 5)(t,x,y)^{<\rho^2>}= \Bbb Q (\sqrt 5)$ $(t,u,v)$.

Step 5.  The action of $\rho$ on $u, \ v$ is given by
$$\rho : u \mapsto \frac {y-\dfrac
by}{\dfrac {ay}x-\dfrac{bx}{y}}, \qquad v \mapsto  \frac {-(x-\dfrac
ax)}{\dfrac {ay}x-\dfrac{bx}{y}}.$$

Define $w=u/(tv)$.  Then $\Bbb Q(\sqrt 5)(t, u, v) = \Bbb Q(\sqrt 5)(t, v,
w)$ and $$\rho: \sqrt 5 \mapsto - \sqrt 5, \ w \mapsto -1/w, \ v \mapsto
\lambda /v$$ where $\lambda = 1/(w-(1/w))$ because
$$\frac {x- \dfrac ax}{\dfrac {ay}x - \dfrac {bx}y} = - \frac
u{bu^2-av^2}. \tag1$$

Note that the above identity is the identity (3) in \cite {CHK, p. 156}.

Define $s= \sqrt 5(1+t)/(1-t)$.  Then $\rho(s)=s$.  Thus $\Bbb Q(\sqrt
5)(t, v, w)^{<\rho>}= \Bbb Q(\sqrt 5)(s, v, w)^{<\rho>}= \Bbb Q(\sqrt
5)(v, w)^{<\rho>}(s)$.  It remains to prove that $\Bbb Q(\sqrt
5)(v,w)^{<\rho>}$ is rational over $\Bbb Q$.

Step 6.  Let $\alpha= \sqrt 5 -2$ and $\beta=1/(w+1)$.  Then $\alpha \cdot
\rho(\alpha)=-1$ and $\beta \cdot \rho(\beta)= 1/(w-(1/w))= \lambda$. Define $W=
w/\alpha, \ V=v/ \beta$.  Then $\Bbb Q(\sqrt 5)(v,w)= \Bbb Q(\sqrt 5)(V,W)$ and
$\rho(V)=1/V$, $\rho(W)=1/W$.  Define $X=(1+V)/(1-V), \ Y=(1+W)/(1-W)$.  Then $\Bbb
Q(\sqrt 5)(V,W)=\Bbb Q(\sqrt 5)(X, Y)$ and $\rho(X)=-X, \ \rho(Y)=-Y$.  Since
$\rho(\sqrt 5X)=\sqrt 5X$ and $\rho(\sqrt 5 Y)=\sqrt 5 Y$, it follows that $\Bbb
Q(\sqrt 5)(X,Y)^{<\rho>}=\Bbb Q(\sqrt 5)(\sqrt 5 X,\sqrt 5Y)^{<\rho>}=\Bbb Q(\sqrt
5)^{<\rho>}(\sqrt 5 X, \sqrt 5 Y)=\Bbb Q(\sqrt 5 X, \sqrt 5 Y)$ is rational over $\Bbb
Q$. \qed \enddemo

\newpage

\centerline {REFERENCES}

\

\roster

\noindent \item "[Br1]" S.  Breuer, {\it Zyklische Gleichungen 6. Grades
und Minimalbasis}, Math. Ann. {\bf 86}(1922) 108-113.

\noindent \item "[Br2]" S.  Breuer, {\it Zur Bestimmung der metazyklischen
 Minimalbasis von Primzahlgrad}, Math. Ann. {\bf 92}(1924) 126-144.

\noindent \item "[Br3]" S.  Breuer, {\it Metazyklischen
 Minimalbasis und komplexe Primzahlen}, J. reine angew. Math {\bf 156}
 (1927) 13-42.

\noindent \item "[CHK]" H.  Chu, S.  J.  Hu and M.  Kang, {\it Noether's
problem for dihedral 2-groups}, Comment. Math. Helv. {\bf 79}(2004)
147-159.

\noindent \item "[DM]"F.\ DeMeyer and T. McKenzie, {\it On generic polynomials}, J.
Algebra {\bf 261}(2003) 327-333.

\noindent \item "[Fi]" E. Fischer, {\it Die Isomorphie der
Invariantenk\"orper der endlichen Abel'schen Gruppen linearer
Transformationen}, Nachr. K\"onigl. Ges. Wiss. G\"ottingen (1915) 77-80.

\noindent \item "[Fu]" P. Furtw\"angler, {\it \"Uber Minimalbasen f\"ur
K\"oper rationalen Funktionen}, S. -B. Akad. Wiss. Wien {\bf 134}(1925)
69-80.

\noindent \item "[Gr]  "
W.\ Gr\"obner,  {\it Minimalbasis der Quaternionengruppe},
Monatshefte f\"ur Math.\ und Physik  {\bf 41}(1934) 78-84.

\noindent \item "[Ha]"
M.\ Hajja  {\it Rationality of finite groups of monomial automorphisms of $K(X,Y)$},
J.\ Algebra {\bf 109}(1987) 46-51.

\noindent \item "[Ka1]"
 M.\ Kang  {\it Introduction to Noether's problem for dihedral groups},
Algebra Colloquium  {\bf 11}(2004) 71-78.

\noindent \item "[Ka2] "
M.\ Kang,  {\it Noether's problem for dihedral 2-groups II},
to appear in "Pacific J. Math.".

\noindent \item "[Ke] "
I. Kersten, {\it Noether's proiblem and normalization},
Jber. Deutsch Math.-Verein.\ {\bf 100} (1999) 3-22.

\noindent \item "[Le]"H. W. Lenstra, Jr., {\it Rational functions
invariant under a finite abelian group},
Invent. math. {\bf 25}(1974) 299-325.

\noindent \item"[Ma] " K. Masuda, {\it Application of the theory of the
group of classes
 of projective modules to the existence problem of independent
 parameters of invariant}, J. Math. Soc. Japan {\bf 20} (1968) 223-232.

\noindent
 \item"[Mi] "  J. Milnor, Introduction to algebraic $K$-theory, Princetion
 Univ. Press, 1971, Princeton.

\noindent \item "[No]" E. Noether, {\it Rationale Funktionenk\"oper},
Jber. Deutsch. Math. -Verein. {\bf 22}(1913) 316-319.

\noindent \item "[Sa1]" D.\ J.\ Saltman, {\it Generic Galois extensions and problems in field theory},
Adv. Math.\ {\bf 43} (1982) 250-283.

\noindent \item "[Sa2]" D.\ J.\ Saltman, {\it Groups acting on fields: Noether's problem},
Contemp. Math.\ {\bf 43} (1985) 267-277.

\noindent \item "[Sw1] " R.\ G.\ Swan, {\it Invariant rational functions and a problem
of Steenrod}, Invent.\ math.\ {\bf 7} (1969) 148-158.

\noindent \item "[Sw2] " R.\ G.\ Swan, {\it Galois Theory}, in ``Emmy Noether, a
tribute to her life and work", Marcel Dekker, 1981, New York.

\noindent \item "[Sw3] " R.\ G.\ Swan, {\it Noether's problem in Galois Theory}, in
``Emmy Noether in Bryn Mawr", Springer-Verlag, 1983, New York.
\endroster

\end